# Counterexample of a theorem in "Wiener index of a fuzzy graph and application to illegal immigration networks"

## Masoud Ghods[1,*], Zahra Rostami [2]

[1]Department of Mathematics, Semnan University, Semnan 35131-19111, Iran, mghods@semnan.ac.ir
[2]Department of Mathematics, Semnan University, Semnan 35131-19111, Iran, zahrarostami.98@semnan.ac.ir

* Correspondence: mghods@semnan.ac.ir; Tel.: (09122310288)

**Abstract**

In the article, we review and critique the Corollary and Theorem of "Wiener index of a fuzzy graph and application to illegal immigration networks", and in addition to providing examples of violations.

**Keywords:** Fuzzy graph; Connectivity index; Wiener index; Fuzzy tree; Fuzzy cycle

1. **Introduction and Preliminaries**

Given that our purpose here is to provide a violation example for a result and a theorem of Article "Wiener index of a fuzzy graph and application to illegal immigration networks". Therefore, the reader is asked to refer to [3].

**Definition 1.1**. [2] Let $G: (\sigma, \mu)$ be a Fuzzy graph. The *Connectivity index* $(CI)$ of G defined by

$$CI(G) = \sum_{u,v \in \sigma^*} \sigma(u)\sigma(v) CONN_G(u,v),$$

where $CONN_G(u,v)$ is the strength of connectedness between $u$ and $v$.

**Definition 1.2**. [1] Let $G: (\sigma, \mu)$ be a Fuzzy graph. The *Wiener index* $(WI)$ of $G$ defined by

$$WI(G) = \sum_{u,v \in \sigma^*} \sigma(u)\sigma(v) d_s(u,v),$$

where $d_s(u,v)$ is the minimum sum of weights of geodesics from $u$ to $v$.

2. **Results**

In this section, we first refer to the Corollary of the article in question and then show with an example that this result is not established.

> **Corollary*** [1] in a fuzzy tree $G: (\sigma, \mu)$ with $F$, $WI(G) = WI(F) = CI(F)$.

**Example 2.1.** (Counterexample) consider the fuzzy graph $G: (\sigma, \mu)$ with $\sigma^* = \{a, b, c, d, e\}$, $\sigma(x) = 1$ for every $x \in \sigma^*$, $\mu(ab) = 0.1, \mu(bc) = \mu(ec) = 0.3, \mu(cd) = 0.5, \mu(ae) = 0.6$. Clearly G is a fuzzy tree and hence there exists the unique maximum spanning tree $(MST)$ for $G$. also there exists unique geodesics between every pair of vertices in G. by calculation

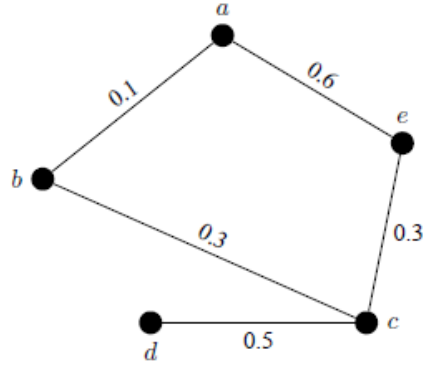

Figure: A fuzzy graph $G: (\sigma, \mu)$ with $\sigma^* = \{a, b, c, d, e\}$

$$WI(G) = d_s(a,b) + d_s(a,c) + d_s(a,d) + d_s(a,e) + d_s(b,c) + d_s(b,d) + d_s(b,e) \\ + d_s(c,d) + d_s(c,e) + d_s(d,e) \\ = 1.2 + 0.9 + 1.4 + 0.6 + 0.3 + 0.8 + 0.6 + 0.5 + 0.3 + 0.8 = 7.4.$$

Clearly $WI(F) = WI(G) = 7.4$. but we have for $CI(F)$
$$CI(F) = CONN_G(a,b) + CONN_G(a,c) + CONN_G(a,d) + CONN_G(a,e) + CONN_G(b,c) \\ + CONN_G(b,d) + CONN_G(b,e) + CONN_G(c,d) + CONN_G(c,e) \\ + CONN_G(d,e) = 0.3 + 0.3 + 0.3 + 0.6 + 0.3 + 0.3 + 0.3 + 0.5 + 0.3 + 0.3 \\ = 3.5.$$

As can be seen $WI(G) = WI(F) \neq CI(F)$.

Now look at the following theorem presented in this article. Then we show with two examples that it is not necessarily established.

**Theorem* [1].** Let $G: (\sigma, \mu)$ be a saturated fuzzy cycle with $C^* = C_n$ of length $n$ such that each $\alpha$-strong edge has strength $\kappa$ and that of each $\beta$-strong edge is $\eta$, then
$$WI(G) = \frac{n[(n+3)^2 - 6]}{16}(\kappa + \eta).$$

☐

**Example 2.2.** Let $G:(\sigma,\mu)$ be a saturated fuzzy cycle with $C^* = C_4$ of length 4 such that each $\alpha$-strong edge has strength $\kappa$ and that of each $\beta$-strong edge is $\eta$, with $\sigma^* = \{a,b,c,d\}$, $\sigma(x) = 1$ for every $x \in \sigma^*$. Let $\mu(ab) = \mu(cd) = \kappa$ and $\mu(bc) = \mu(ad) = \eta$. Then

$$WI(G) = d_s(a,b) + d_s(a,c) + d_s(a,d) + d_s(b,c) + d_s(b,d) + d_s(c,d)$$
$$= \kappa + (\kappa + \eta) + \eta + \eta + (\eta + \kappa) + \kappa = 4\kappa + 4\eta = 4(\kappa + \eta).$$

But, by use *Theorem\**,

$$WI(G) = \frac{n[(n+3)^2 - 6]}{16}(\kappa + \eta) = \frac{4[(4+3)^2 - 6]}{16}(\kappa + \eta) = \frac{4[49 - 6]}{16}(\kappa + \eta)$$
$$= \frac{43}{4}(\kappa + \eta).$$

As can be seen, the value obtained is not correct.

Now, in the example blow let n=6.

**Example 2.3.** Let $G:(\sigma,\mu)$ be a saturated fuzzy cycle with $C^* = C_6$ of length 4 such that each $\alpha$-strong edge has strength $\kappa$ and that of each $\beta$-strong edge is $\eta$, with $\sigma^* = \{a,b,c,d,e,f\}$, $\sigma(x) = 1$ for every $x \in \sigma^*$. Let $\mu(ab) = \mu(cd) = \mu(ef) = \kappa$ and $\mu(bc) = \mu(de) = \mu(af) = \eta$.

|   | a | b | c | d | e | f |
|---|---|---|---|---|---|---|
| a | - | $\kappa$ | $\kappa + \eta$ | $2\eta + \kappa$ | $\eta + \kappa$ | $\eta$ |
| b | $\kappa$ | - | $\eta$ | $\eta + \kappa$ | $2\eta + \kappa$ | $\kappa + \eta$ |
| c | $\kappa + \eta$ | $\eta$ | - | $\kappa$ | $\kappa + \eta$ | $2\eta + \kappa$ |
| d | $2\eta + \kappa$ | $\eta + \kappa$ | $\kappa$ | - | $\eta$ | $\eta + \kappa$ |
| e | $\eta + \kappa$ | $2\eta + \kappa$ | $\kappa + \eta$ | $\eta$ | - | $\kappa$ |
| f | $\eta$ | $\kappa + \eta$ | $2\eta + \kappa$ | $\eta + \kappa$ | $\kappa$ | - |

Then

$$WI(G) = d_s(a,b) + d_s(a,c) + d_s(a,d) + d_s(a,e) + d_s(a,f) + d_s(b,c) + d_s(b,d) + d_s(b,e)$$
$$+ d_s(b,f) + d_s(c,d) + d_s(c,e) + d_s(c,f) + d_s(d,e) + d_s(d,f) + d_s(e,f)$$
$$= \kappa + (\kappa + \eta) + (2\eta + \kappa) + (\eta + \kappa) + (\eta) + (\eta) + (\eta + \kappa) + (2\eta + \kappa) + (\kappa + \eta)$$
$$+ (\kappa) + (\kappa + \eta) + (2\eta + \kappa) + (\eta) + (\eta + \kappa) + (\kappa) = 12\kappa + 15\eta.$$

But, by use *Theorem\**, we have

$$WI(G) = \frac{n[(n+3)^2 - 6]}{16}(\kappa + \eta) = \frac{6[(6+3)^2 - 6]}{16}(\kappa + \eta) = \frac{6[81 - 6]}{16}(\kappa + \eta)$$
$$= \frac{6 \times 75}{16}(\kappa + \eta).$$

As can be seen, the value obtained is not correct.

**Refrences**